\documentclass[final]{siamltex}
\usepackage{graphicx}

\newtheorem{pbm}[theorem]{Open Problem}
\newtheorem{rem}[theorem]{Remark}
\newcommand{\conv}{{\it conv}}
\newcommand{\link}{{\it link}}
\renewcommand{\max}{{\it max}}
\renewcommand{\min}{{\it min}}
\newcommand{\vertices}{{\it vertices}}

\def\bbbr{\mbox{\boldmath$R$}}
\title{Extremal properties for dissections of\\convex 3-polytopes}

\author{%
Jes\'us A. De Loera%
\thanks{Dept. of Mathematics, Univ. of California-Davis
({\tt deloera@math.ucdavis.edu}). The research of this author
partially supported by NSF grant DMS-0073815.}
\and
Francisco Santos%
\thanks{Depto. de Matem\'aticas, Estad. y Comput., 
Univ. de Cantabria
({\tt santos@matesco.unican.es}).
The research of this author was supported partially by grant PB97-0358
of the Spanish Direcci\'on General de Investigaci\'on
Cient\'{\i}fica y T\'ecnica.}
\and
Fumihiko Takeuchi%
\thanks{Dept. of Information Science, Univ. of Tokyo
({\tt fumi@is.s.u-tokyo.ac.jp}).}}

\begin{document}

\maketitle

\begin{abstract}
A dissection of a convex $d$-polytope is a partition of the polytope
into $d$-simplices whose vertices are among the vertices of the
polytope. Triangulations are dissections that have the additional
property that the set of all its simplices forms a simplicial
complex. The size of a dissection is the number of $d$-simplices it
contains. This paper compares triangulations of maximal size with
dissections of maximal size. We also exhibit lower and upper bounds
for the size of dissections of a $3$-polytope and analyze extremal
size triangulations for specific non-simplicial polytopes: prisms,
antiprisms, Archimedean solids, and combinatorial $d$-cubes.
\end{abstract}



\pagestyle{myheadings}
\thispagestyle{plain}

\section{Introduction}

Let ${\cal A}$ be a point configuration in $\bbbr^d$ with its convex
hull $\conv({\cal A})$ having dimension $d$. A set of $d$-simplices
with vertices in ${\cal A}$ is a {\em dissection} of ${\cal A}$ if no
pair of simplices has an interior point in common and their union
equals $\conv({\cal A})$.  A dissection is a {\em triangulation } of
${\cal A}$ if in addition any pair of simplices intersects at a common
face (possibly empty). The {\em size} of a dissection is the number of
$d$-simplices it contains. We say that a dissection is {\em
mismatching} when it is not a triangulation (i.e.\ it is not a
simplicial complex). In this paper we study mismatching dissections of
maximal possible size for a convex polytope and compare them with
maximal triangulations. This investigation is related to the study of
Hilbert bases and the hierarchy of covering properties for polyhedral
cones which is relevant in Algebraic Geometry and Integer Programming
(see \cite{bruns,firzie,sebo}). Maximal dissections are relevant
also in the enumeration of interior lattice points and its
applications (see \cite{barvinok,kantor} and references there).

It was first shown by Lagarias and Ziegler that dissections of maximal
size turn out to be, in general, larger than maximal triangulations,
but their example uses interior points \cite{lagzie}. Similar
investigations were undertaken for mismatching minimal dissections and
minimal triangulations of convex polytopes \cite{brehmo}. In this
paper we augment previous results by showing that it is possible to
have {\em simultaneously}, in the same $3$-polytope, that the size of
a mismatching minimal (maximal) dissection is smaller (larger) than
any minimal (maximal) triangulation.  In addition, we show that the
gap between the size of a mismatching maximal dissection and a maximal
triangulation can grow linearly on the number of vertices and that
this occurs already for a family of simplicial convex $3$-polytopes.
A natural question is how different are the upper and lower bounds for
the size of mismatching dissections versus those bounds known for
triangulations (see \cite{rotstraus}).  We prove lower and upper
bounds on their size with respect to the number of vertices for
dimension three and exhibit examples showing that our technique of
proof fails already in dimension four. Here is the first summary of
results:

\begin{theorem}
\label{main}
\label{boundsul}
\begin{enumerate}
  
\item There exists an infinite family of convex simplicial
  $3$-polytopes with increasing number of vertices whose mismatching
  maximal dissections are larger than their maximal triangulations.
  This gap is linear in the number of vertices
  $($Corollary \ref{cor.Pm}\/$)$.

\item \begin{enumerate}

  \item There exists a lattice $3$-polytope with $8$ vertices
    containing no other lattice point other than its vertices whose
    maximal dissection is larger than its maximal triangulations.

  \item There exists a $3$-polytope with $8$ vertices for which,
    simultaneously, its minimal dissection is smaller than minimal
    triangulations and maximal dissection is larger than maximal
    triangulations.

\end{enumerate}
$($Proposition \ref{prop.latticepoly}\/$)$

\item If $D$ is a mismatching dissection of a $3$-polytope with $n$
  vertices, then the size of $D$ is at least $n-2$.
  In addition, the size of $D$ is bounded above by
${n-2 \choose 2}$
  $($Proposition \ref{bounds}\/$)$.
\end{enumerate}
\end{theorem}

A consequence of our third point is that the result of
\cite{brehmo}, stating a linear gap between the size of minimal
dissections and minimal triangulations, is best possible.  The results
are discussed in Sections \ref{maxdissvsmaxt} and \ref{upperlower}.

The last section presents a study of maximal and minimal
triangulations for combinatorial $d$-cubes, three-dimensional prisms
and anti-prisms, as well as other Archimedean polytopes. The following
theorem and table summarize the main results:

\begin{theorem}
\begin{enumerate}
\item There is a constant $c>1$ such that for every $d\ge3$
the maximal triangulation among all possible combinatorial $d$-cubes
has size at least $c^d d!$
$($Proposition \ref{combinatorialcube}\/$)$.

\item For a three-dimensional $m$-prism, in any of its possible
  coordinatizations, the size of a minimal triangulation is
  $2m-5+\lceil {m\over 2}\rceil$. For an $m$-antiprism, in any of its
  possible coordinatizations, the size of a minimal triangulation is
  $3m-5$ $($Proposition \ref{prismantip}\/$)$. The size of a maximal
  triangulation of an $m$-prism depends on the coordinatization, and in
  certain natural cases it is $(m^2+m-6)/2$ $($Proposition
  \ref{maxprism}\/$)$.

\item The following table specifies sizes of the minimal and maximal
triangulations for some Platonic and Archimidean solids. These results were
obtained via integer programming calculations using the approach
described in \cite{dhss}. All computations used the canonical
symmetric coordinatizations for these polytopes \cite{coxeter}.  
The number of vertices is indicated in
parenthesis
$($Remark \ref{rem.dhss}\/$)$:

\begin{table}[h]
\begin{center}
\footnotesize
\font\ninerm=cmr8
\ninerm
\begin{tabular}{|c|c|c|}
\hline
 {\bf $P$}   & $|T_{\min}(P)|$ & $|T_{\max}(P)|$ \cr
\hline
 Icosahedron (12)               & 15 &  20  \cr
 Dodecahedron (20)              & 23 &  36  \cr
 Cuboctahedron (12)             & 13 &  17  \cr
 Icosidodecahedron (30)         & 45 & \normalsize ?  \cr
 Truncated Tetrahedron (12)     & 10 &  13  \cr
 Truncated Octahedron (24)      & 27 & \normalsize ?  \cr
 Truncated Cube (24)            & 25 &  48  \cr
 Small Rhombicuboctahedron (24) & 35 & \normalsize ?  \cr
 Pentakis Dodecahedron (32)     & 54 & \normalsize ?  \cr
 Rhombododecahedron (14)        & 12 &  21  \cr
\hline
\end{tabular}
\end{center}
\caption{Sizes of extremal triangulations of Platonic and Archimidean solids.}
\label{res_table}
\end{table}

\end{enumerate}

\end{theorem}
\section{Maximal dissections of $\mathbf 3$-polytopes} \label{maxdissvsmaxt}

\begin{figure}[htb]
  \begin{center}
        \includegraphics[scale=.4]{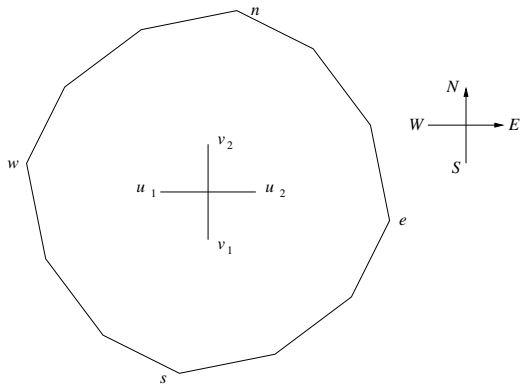}
  \end{center}
\caption{North, South, East, and West vertices.}
 \label{fig.NSEW}
\end{figure}

We introduce some important definitions and conventions: We denote by
$Q_m$ a convex $m$-gon with $m$ an even positive integer. Let
$v_1v_2$ and $u_1u_2$ be two edges parallel to $Q_m$, orthogonal to
each other, on opposite sides of the plane containing $Q_m$, and such
that the four segments $v_iu_j$ intersect the interior of $Q_m$. We
suppose that $v_1v_2$ and $u_1u_2$ are not parallel to any diagonal or edge
of $Q_m$.  The convex hull $P_m$ of these points has $m+4$
vertices and it is a simplicial polytope.  We will call north (respectively
south) vertex of $Q_m$ the one which maximizes (respectively minimizes)
the scalar product with the vector $v_2-v_1$. Similarly, we will call
east (west) the vertex which maximizes (minimizes) the scalar product
with $u_2-u_1$. We denote these four vertices $n$, $s$, $e$ and $w$,
respectively. See Figure \ref{fig.NSEW}.

\begin{figure}[htb]
  \begin{center}
        \includegraphics[scale=.4]{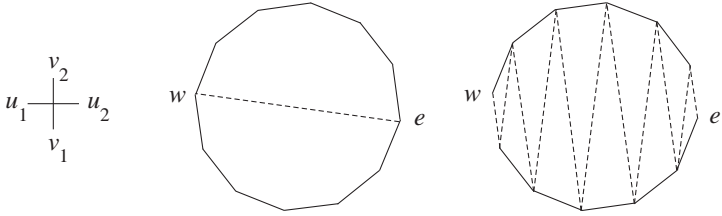}        
  \end{center}
\caption{The minimal monotone path $($middle$)$ and
the maximal monotone path made by the vertical edges $($right$)$
in the direction $u_1u_2$.} \label{zigzag}
\end{figure}

We say that a directed path of edges inside $Q_m$ is {\em monotone in
the direction} $v_1v_2$ (respectively $u_1u_2$) when the vertices of
the path appear in the path following the same order given by the
scalar product with $v_2-v_1$ (respectively $u_2-u_1$). An equivalent
formulation is that any line orthogonal to $v_1v_2$ cuts the path in at
most one point. We remark that by our choice of $v_1v_2$ and $u_1u_2$
all vertices of $Q_m$ are ordered by the values of their scalar products
with $v_2-v_1$ and also with respect to $u_2-u_1$.  In the same way, a
sequence of vertices of $Q_m$ is {\em ordered in the direction} of
$v_1v_2$ (respectively $u_1u_2$), if the order is the same as the one
provided by using the values of the scalar products of the points with
the vector $v_2-v_1$ (respectively $u_2-u_1$).  Consider the two
orderings induced by the directions of $v_1v_2$ and $u_1u_2$ on the
set of vertices of $Q_m$. Let us call {\em horizontal} (respectively {\em
vertical}\/) any edge joining two consecutive vertices in the
direction of $v_1v_2$ (respectively of $u_1u_2$). As an example, 
if $Q_m$ is regular then
the vertical edges in $Q_m$ form a zig-zag path as shown 
in Figure \ref{zigzag}.

\begin{figure}[htb]
  \begin{center}
        \includegraphics[scale=.35]{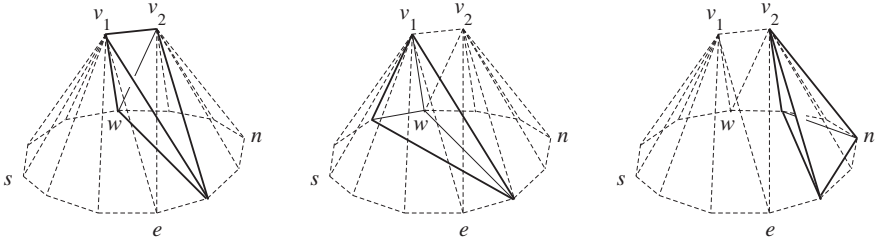}
  \end{center}
\caption{Three types of tetrahedra in $R_m$.}
 \label{fig.Rm}
\end{figure}

Our examples in this section will be based on the following
observation and are inspired by a similar analysis of maximal dissections
of dilated empty lattice tetrahedra in $\bbbr^3$ by Lagarias
and Ziegler \cite{lagzie}: Let $R_m$ be the convex hull of the $m+2$ vertices
consisting of the $m$-gon $Q_m$ and $v_1,v_2$. $R_m$ is exactly one
half of the polytope $P_m$.  Consider a triangulation $T_0$ of $Q_m$
and a path $\Gamma$ of edges of $T_0$ monotone with respect to the
direction $u_1u_2$. Observe that $\Gamma$
divides $T_0$ in two regions, which we will call the ``north'' and the
``south''.
Then, the following
three families of tetrahedra form a triangulation $T$ of $R_m$: the
edges of $\Gamma$ joined to the edge $v_1v_2$; 
the southern triangles of $T_0$ joined to $v_1$;
and the northern triangles of $T_0$ joined to $v_2$
(see Figure \ref{fig.Rm}).
Moreover, all the triangulations of $R_m$ are obtained in this
way: Any triangulation $T$ of $R_m$ induces a triangulation $T_0$ of
$Q_m$. The link of $v_1v_2$ in $T$ is a monotone path of edges contained
in $T_0$ and it divides $T_0$ in two regions, joined respectively to $v_1$
and $v_2$.

Using the Cayley trick, one can also think of the
triangulations of $R_m$ as the fine mixed subdivisions of the
Minkowski sum $Q_m + v_1v_2$ (see \cite{cayley} and references
within).

The size of a triangulation of $R_m$ equals $m-2+\vert \Gamma\vert$,
where $\vert \Gamma\vert$ is the number of edges in the path
$\Gamma$. There is a unique minimal path in $Q_m$ of
length one (Figure \ref{zigzag}, middle)
and a unique maximal path of length $m-1$ (Figure \ref{zigzag},
right). 
Hence the
minimal and maximal triangulations of $R_m$ have, respectively, $m-1$
and $2m-3$ tetrahedra.  The maximal triangulation is unique, but the
minimal one is not: after choosing the diagonal in $\Gamma$ the rest
of the polygon $Q_m$ can be triangulated in many ways.  From the above
discussion regarding $R_m$ we see that we could independently
triangulate each of the two halves of $P_m$ with any number of
tetrahedra from $m-1$ to $2m-3$. Hence, $P_m$ has dissections of sizes
going from $2m-2$ to $4m-6$.  Among the triangulations of $P_m$, we
will call {\em halving triangulations} those that triangulate the two
halves of $P_m$. Equivalently, the halving triangulations are those
which do not contain any of the four edges $v_iu_j$.

\begin{proposition}
\label{prop.nonhalving}
Let $P_m$ be as described above, with $Q_m$ being a regular $m$-gon.
No triangulation of $P_m$ has
more than ${7m\over 2}+1$ tetrahedra. On the other hand, there are
mismatching dissections of $P_m$ with $4m-6$ tetrahedra.
\end{proposition}
\begin{proof}
  Let $T$ be a triangulation of $P_m$.  It is an easy application of
  Euler's formulas for the 3-ball and 2-sphere that the number of
  tetrahedra in a triangulation of any 3-ball without interior
  vertices equals the number of vertices plus interior edges minus
  three (such formula appears for instance in \cite{Edels}). Hence our
  task is to prove that $T$ has at most ${5 m\over 2}$ interior edges.
  For this, we classify the interior edges according to how many
  vertices of $Q_m$ they are incident to. There are only four edges
  not incident to any vertex of $Q_m$ (the edges $v_iu_j$,
  $i,j\in\{1,2\}$). Moreover, $T$ contains at most $m-3$ edges
  incident to two vertices of $Q_m$ (i.e.\ diagonals of $Q_m$), since
  in any family of more than $m-3$ such edges there are pairs which
  cross each other. Thus, it suffices to prove that $T$ contains at most
  ${3 m\over 2}-1$ edges incident to just one vertex of $Q_m$, i.e.\ of
  the form $v_ip$ or $u_ip$ with $p\in Q_m$.

Let $p$ be any vertex of $Q_m$. If $p$ equals $w$ or $e$ then the
edges $pv_1$ and $pv_2$ are both in the boundary of $P_m$; for any
other $p$, exactly one of $pv_1$ and $pv_2$ is on the boundary and the
other one is interior. Moreover, we claim that if $pv_i$ is an
interior edge in a triangulation $T$, then the triangle $pv_1v_2$
appears in $T$. This is so because there is a plane containing $pv_i$
and having $v_{3-i}$ as the unique vertex on one side. At the same
time the link of $pv_i$ is a cycle going around the edge. Hence,
$v_{3-i}$ must appear in the link of $pv_i$.  It follows from the
above claim that the
number of interior edges of the form $pv_i$ in $T$ equals the number
of vertices of $Q_m$ other than $w$ and $e$ in the link of $v_1v_2$.
In a similar way, the number of interior edges of the form $pu_i$ in
$T$ equals the number of vertices of $Q_m$ other than $n$ and $s$ in
the link of $u_1u_2$. In other words, if we call
$\Gamma_u=\link_T(v_1v_2)\cap Q_m$ and $\Gamma_v=\link_T(u_1u_2)\cap
Q_m$ (the $u$, $v$ in the index and of the vertices are reversed,
because in this way $\Gamma_u$ is monotone with respect to $u_1u_2$,
and $\Gamma_v$ with respect to $v_1v_2$), then
the number of interior edges in $T$ incident to exactly one vertex of
$Q_m$ equals $|\vertices(\Gamma_v)|+|\vertices(\Gamma_u)|-4$. Our goal
is to bound this number. 
As an example, Figure \ref{fig.gamma} shows the intersection of
$Q_m$ with a certain triangulation of $P_m$ ($m=12$). The link of $v_1v_2$
in this triangulation is the chain of vertices and edges
$wabu_1nu_2ce$ (the star of $v_1v_2$ is 
marked in thick and grey in the figure). 
$\Gamma_u$ consists of the chains $wab$ and $ce$ and
the isolated vertex $n$. In turn, the link of $u_1u_2$ is the chain
$nv_1s$ and $\Gamma_v$ consists of the isolated vertices $n$ and $s$.


\begin{figure}[htb]
  \begin{center}
        \includegraphics[scale=.4]{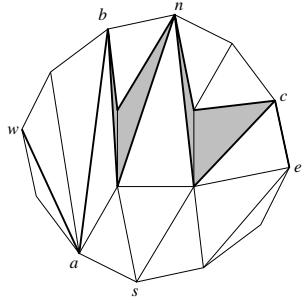}
  \end{center}
\caption{Illustration of the proof of Proposition \ref{prop.nonhalving}.}
 \label{fig.gamma}
\end{figure}

Observe that $\Gamma_v$ has at most three connected components,
because it is obtained by removing from $\link_T(u_1u_2)$ (a path) the
parts of it incident to $v_1$ and $v_2$, if any.
Each component is monotone in the direction of
$v_1v_2$ and the projections of any two components to a line parallel
to $v_1v_2$ do not overlap.  The sequence of vertices of $Q_m$ ordered
in the direction of $v_1v_2$, can have a pair of consecutive vertices
contained in $\Gamma_v$ only where there is a horizontal edge in
$\Gamma_v$ or in the at most two discontinuities of $\Gamma_v$. This
is true because $Q_m$ is a regular $m$-gon.

We denote $n_{hor}$ the number of horizontal edges in $\Gamma_v$ and $n_{hor}'$
this number plus the number of discontinuities in $\Gamma_v$ (hence
$n_{hor}'\le n_{hor} + 2$). Every non-horizontal edge of $\Gamma_v$ produces a
jump of at least two in the $v_1v_2$-ordering of the vertices of
$P_m$, hence we have

$$\vert \vertices(\Gamma_v)\vert -1 -n_{hor}' \le {m-1-n_{hor}'\over 2}.$$

Analogously, and with the obvious similar meaning for $n_{vert}$ and $n_{vert}'$,

$$
\vert \vertices(\Gamma_u)\vert -1 -n_{vert}' \le {m-1-n_{vert}'\over 2}.
$$

Since $\Gamma_u\cup \Gamma_v$ can be completed to a triangulation of
$Q_m$, and exactly four non-interior edges of $Q_m$ are horizontal or
vertical, we have $n_{hor}+n_{vert}\le (m-3)+4 = m+1$, i.e.\
$n_{hor}'+n_{vert}'\le m+5$. Hence,

$$
\vert \vertices(\Gamma_v)\vert + \vert \vertices(\Gamma_u)\vert
\le \left\lfloor{2m+2+n_{hor}'+n_{vert}'\over 2}\right\rfloor
\le \left\lfloor{3m+7\over 2}\right\rfloor = {3m\over 2} +3.
$$

Thus, there are at most ${3m\over 2}-1$ interior edges in $T$ of the
form $pv_i$ or $pu_i$ and at most ${5m\over 2}$ interior edges in
total, as desired.  
\end{proof}

\begin{corollary}
\label{cor.Pm}
The polytope $P_m$ described above has the following properties:
\begin{itemize}
\item It is a simplicial $3$-polytope with $m+4$ vertices.
\item Its maximal dissection has at least $4m-6$ tetrahedra.
\item Its maximal triangulation has at most ${7m\over 2}+1$ tetrahedra.
\end{itemize}
In particular, the gap between sizes of the maximal dissection and
maximal triangulation is linear on the number of vertices.
\end{corollary}

\begin{figure}[htb]
  \begin{center}
\includegraphics[scale=.25]{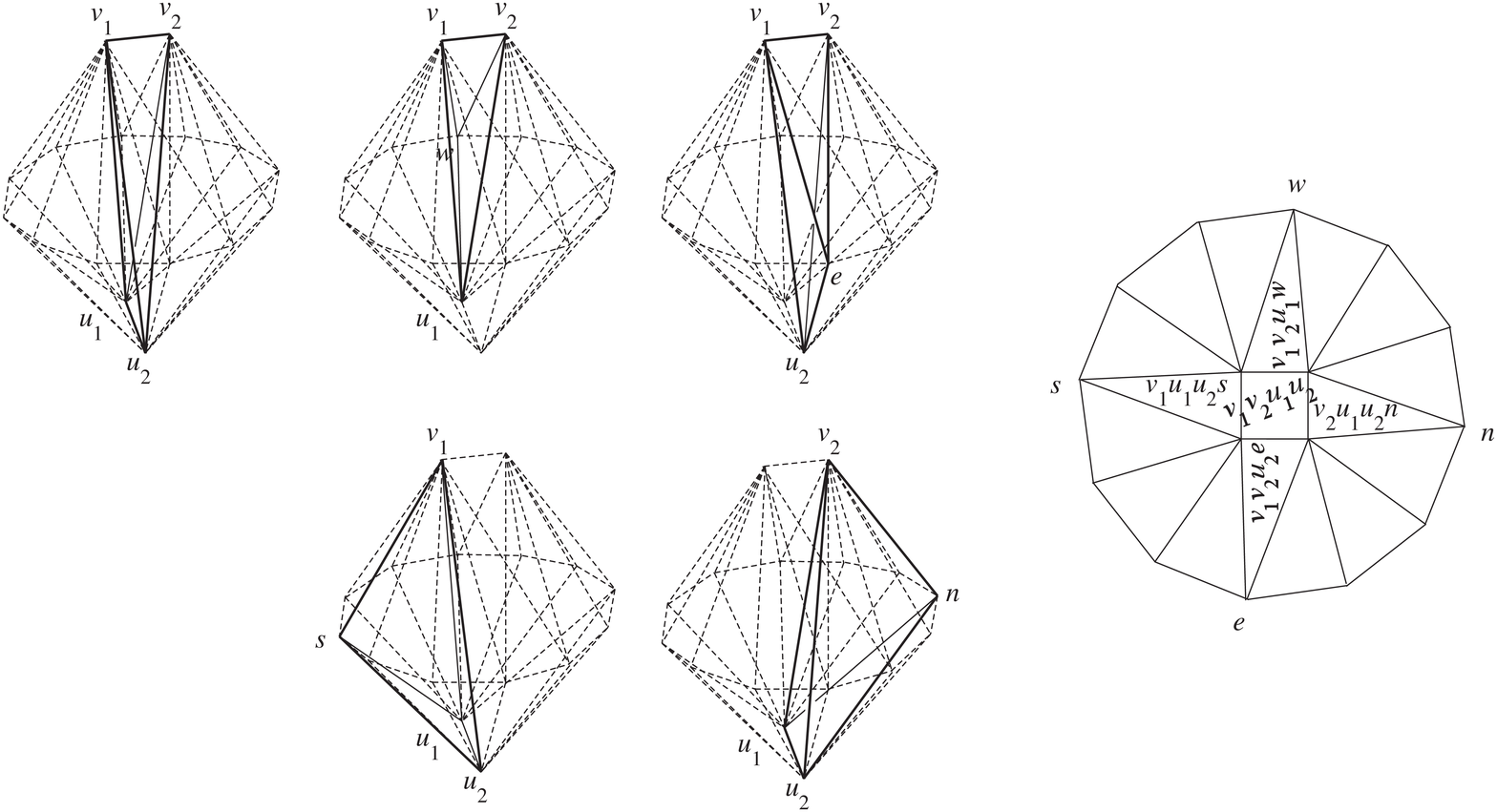}               
  \end{center}
\caption{For the triangulation of $P_m$ with $m+5$ tetrahedra,
  its five central tetrahedra $($left$)$ and the intersection 
  of the triangulation with the polygon $Q_m$ $($right$)$ are shown.
  The four interior vertices are the intersection points of the edges
  $v_1u_1$, $v_1u_2$, $v_2u_1$ and $v_2u_2$ with the plane containing $Q_m$.}
 \label{fig.Pmsmalltri}
\end{figure}

Three remarks are in order: First, the size of the maximal triangulation
for $P_m$ may depend on the coordinates or, more specifically on which diagonals
of $Q_m$ intersect the tetrahedron $v_1v_2u_1u_2$. Second, concerning
the size of the minimal triangulation of $P_m$, we can easily describe
a triangulation of $P_m$ with only $m+5$ tetrahedra: let the vertices
$n$, $s$, $e$ and $w$ be as defined above  (see Figure
\ref{fig.NSEW}) and let us call northeast, northwest, southeast and
southwest the edges in the arcs $ne$, $nw$, $se$ and $sw$ in the
boundary of $Q_m$. Then, the triangulation consists of the five
tetrahedra $v_1v_2u_1u_2$, $v_1v_2u_1w$, $v_1v_2u_2e$, $v_1u_1u_2s$
and $v_2u_1u_2n$ (shown in the left part of Figure \ref{fig.Pmsmalltri})
together with the edges $v_2u_2$, $v_2u_1$, $v_1u_2$
and $v_1u_1$ joined, respectively, to the northeast, northwest,
southeast and southwest edges of $Q_m$. The right part of Figure
\ref{fig.Pmsmalltri} shows the result of slicing through the
triangulation by the plane containing the polygon $Q_m$.

Finally, although the corollary above states a difference between
maximal dissections and maximal triangulations only for $P_m$ with
$m>14$, experimentally we have observed there is a gap already for
$m=8$. Now we discuss two other interesting examples. The following
proposition constitutes the proof of Theorem \ref{main} (2). 

\begin{proposition} \label{prop.latticepoly}
\begin{enumerate}

\item Consider the following eight points in $\bbbr^3$:
\begin{itemize}
\item The vertices $s=(0,0,0)$, $e=(1,0,0)$, $w=(0,1,0)$ and
$n=(1,1,0)$ of  a square in the plane $z=0$.
\item The vertices $v_1=(-1,0,1)$ and $v_2=(1,1,1)$ 
of a horizontal edge above the square, and
\item The vertices $u_1=(0,1,-1)$ and $u_2=(2,0,-1)$ 
of a horizontal edge below the square.
\end{itemize}
These eight points are the vertices
of a polytope $P$ whose only integer points are precisely its
eight vertices and with the following properties:
\begin{enumerate}
\item Its $($unique\/$)$ maximal dissection has $12$ tetrahedra. All of
them are unimodular, i.e.\ they have volume $1/6$.
\item Its $($several\/$)$ maximal triangulations have $11$ tetrahedra.
\end{enumerate}

\item 
  For the $3$-polytope with vertices 
  $u_1=(1,0,0)$, $w=(1,0,1)$, $v_1=(-1,0,0)$, $s=(-1,0,-1)$, 
  $v_2=(0,1,1)$, $n=(1,1,1)$, $u_2=(0,1,-1)$, $e=(-1,1,-1)$,
  the sizes of its $($unique\/$)$ minimal
  dissection and $($several\/$)$ minimal triangulations
  are $6$ and $7$ respectively, and the sizes of its $($several\/$)$ maximal
  triangulations and $($unique\/$)$ maximal dissection
  are $9$ and $10$ respectively.  
\end{enumerate}
\end{proposition}

\begin{proof} The polytopes
constructed are quite similar to $P_4$ constructed earlier except that
$Q_4$ is non-regular (in part 2) and the segments $u_1u_2$ and
$v_1v_2$ are longer and are not orthogonal, thus ending with different
polytopes.
The polytopes are shown in Figure \ref{fig.prop23(1)}. Figure 
\ref{fig.prop23} describes a maximal dissection of each of them, in
five parallel slices.
Observe that both polytopes have four vertices in
the plane $y=0$ and another four in the plane $y=1$. Hence, the first
and last slices in parts (a) and (b) of Figure \ref{fig.prop23}
completely describe the polytope.
\begin{figure}[htb]
  \begin{center}
        \includegraphics[scale=.52]{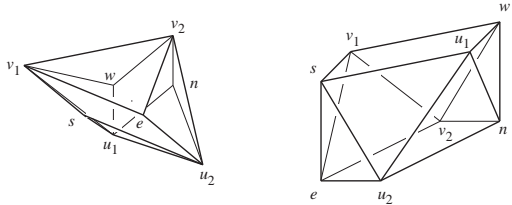}
  \end{center}
\caption{The two polytopes in Proposition \ref{prop.latticepoly}.}
\label{fig.prop23(1)}
\end{figure}

\begin{figure}[htb]
  \begin{center}
        \includegraphics[scale=.8]{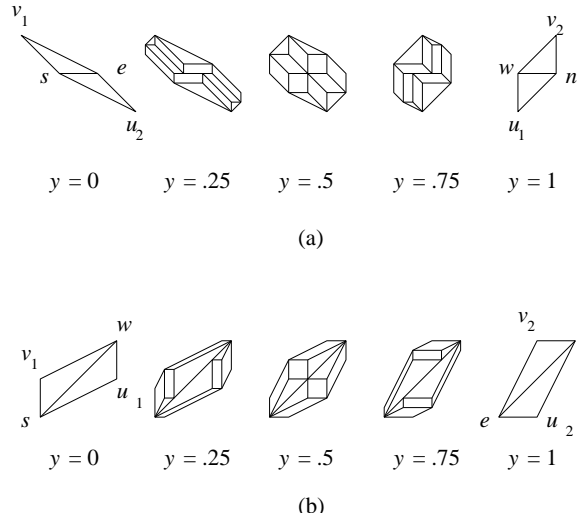}
  \end{center}
\caption{
Five $2$-dimensional slices of the maximal dissections of the polytopes
in Proposition \ref{prop.latticepoly}.
The first and last slices are two facets of the
polytopes containing all the vertices.}  \label{fig.prop23}
\end{figure}

(1) The vertices in the planes $y=0$ and $y=1$ 
form convex quadrangles whose
only integer points are the four vertices. This proves that the eight
points are in convex position and that the polytope $P$ contains no
integer point other than its vertices.  Let us now prove the
assertions on maximal dissections and triangulations of $P$:

(a) Consider the paths of length three $\Gamma_v=\{esnw\}$ and
$\Gamma_u=\{sewn\}$, which are monotone respectively in the directions
orthogonal to $v_1v_2$ and $u_1u_2$. Using them, we can construct two
triangulations of size five of the polytopes $\conv(nsewv_1v_2)$ and
$\conv(nsewu_1u_2)$, respectively. But they do not fill $P$
completely. There is space left for the tetrahedra $swv_1u_1$ and
$env_2u_2$. This gives a dissection of $P$ with twelve tetrahedra. All
the tetrahedra are unimodular, so no bigger dissection is possible.

(b) A triangulation of size 11 can be obtained using the same idea as
above, but with paths $\Gamma_v$ and $\Gamma_u$ of lengths three and
two respectively, which can be taken from the same triangulation of
the square $nswe$.  

To prove that no triangulation has bigger size, it
suffices to show that $P$ does not have any unimodular triangulation.
This means all
tetrahedra have volume $1/6$. We start by recalling a well-known fact
(see Corollary 4.5 in \cite{sebo2}). A lattice tetrahedron has volume
$1/6$ if and only if each of its vertices $v$ lies in a {\em consecutive}
lattice plane parallel to the supporting plane of the opposite facet
to $v$. Two parallel planes are said to be consecutive if their
equations are $ax+by+cz=d$ and $ax+by+cz=d-1$.

Suppose that $T$ is a unimodular triangulation of $P$. 
We will first prove that the triangle
$u_1u_2e$ is in $T$.  The triangular facet $u_1u_2s$ of $P$, lying in
the hyperplane $x+2y+2z=0$, has to be joined to a vertex in the plane
$x+2y+2z=1$. The two possibilities are $e$ and $v_1$. With the same
argument, if the tetrahedron $u_1u_2sv_1$ is in $T$, its facet
$u_1u_2v_1$, which lies in the hyperplane $2x+4y+3z=1$, will be joined
to a vertex in $2x+4y+3z=2$, and the only one is $e$. This finishes
the proof that $u_1u_2e$ is a triangle in $T$.  Now, $u_1u_2e$ is in
the plane $x+2y+z=1$ and must be joined to a vertex in $x+2y+z=2$,
i.e.\ to $w$. Hence $u_1u_2ew$ is in $T$ and, in particular, $T$ uses
the edge $ew$.  $P$ is symmetric under the rotation of order two on
the axis $\{z=0,x={1\over 2}\}$. Applying this symmetry to the
previous arguments we conclude that $T$ uses the edge $ns$ too. But
this is impossible since the edges $ns$ and $ew$ cross each other.

(2) This polytope almost fits the description of $P_4$, except for the
fact that the edges $v_1u_1, v_2u_2$ intersect the boundary and not the
interior of the planar quadrangle $nsew$. With the general techniques
we have described, it is easy to construct halving dissections of this
polytope with sizes from 6 to 10. Combinatorially, the polytope is a
4-antiprism. Hence, Proposition \ref{prismantip} shows that its
minimal triangulation has 7 tetrahedra. The rest of the assertions in
the statement were proved using the integer programming approach
proposed in \cite{dhss}, which we describe in Remark
\ref{rem.dhss}. We have also verified them by enumerating all
triangulations \cite{rambau,imaio}. It is interesting to observe that
if we perturb the coordinates a little so that the planar
quadrilateral $u_1v_1u_2e$ becomes a tetrahedron with the right
orientation and without changing the face lattice of the polytope,
then the following becomes a triangulation with ten tetrahedra:
$\{u_1u_2se,$ $u_1u_2ev_1,$ $u_1u_2v_1w,$ $u_1u_2wn,$ $v_1v_2en,$
$v_1v_2nw,$ $u_1v_1se,$ $v_1u_2ew,$ $u_2wne,$ $v_1wne\}$.
\end{proof}

\section{Bounds for the size of a dissection} \label{upperlower}

Let $D$ be a dissection of a $d$-polytope $P$.  Say two
$(d-1)$-simplices $S_1$ and $S_2$ of $D$ {\em intersect improperly} in
a $(d-1)$-hyperplane $H$ if both lie in $H$, are not identical, and
they intersect with non-empty relative interior.  Consider the following
auxiliary graph: take as nodes the $(d-1)$-simplices of a dissection,
and say that two $(d-1)$-simplices are adjacent if they intersect improperly
in certain hyperplane.  A {\em mismatched region} is the subset of
$\bbbr^d$ that is the union of $(d-1)$-simplices over a connected
component of size larger than one in such a graph.  Later, in
Proposition \ref{badb} we will show some of the complications that can
occur in higher dimensions.

Define the {\em simplicial complex of a dissection} as
all the simplices of the dissection together with their faces, where only faces
that are identical (in $\bbbr^d$) are identified. This construction 
corresponds intuitively to an {\em inflation\/} of the dissection where for 
each mismatched region we move the two groups of $(d-1)$-simplices 
slightly apart leaving the relative boundary of the mismatched region 
joined. Clearly, the simplicial complex of a dissection may be not
homeomorphic to a ball.

The deformed $d$-simplices intersect properly, and the mismatched
regions become holes.  The numbers of vertices and $d$-simplices do
not change. 

\begin{lemma} All mismatched regions for a dissection of  a
convex $3$-polytope $P$ are convex polygons with all vertices among
the vertices of $P$. Distinct mismatched regions have disjoint
relative interiors.
\label{2dmismatch} 
\end{lemma}
\begin{proof}
Let $Q$ be a mismatched region and $H$ the plane containing it.  Since
a mismatched region is a union of overlapping triangles, it is a
polygon in $H$ with a connected interior. If two triangles forming
the mismatched region have interior points in common, they should be
facets of tetrahedra in different sides of $H$.  Otherwise, the two
tetrahedra would have interior points in common, contradicting the
definition of dissection.  Triangles which are facets of tetrahedra in
one side of $H$ cover $Q$. Triangles coming from the other side of $H$
also cover $Q$.

Now, take triangles coming from one side.  As mentioned above, they
have no interior points in common.  Their vertices are among the
vertices of the tetrahedra in the dissection, thus among the vertices
of the polytope $P$.  Hence, the vertices of the triangles are in
convex position, thus the triangles are forming a triangulation of a
convex polygon in $H$ whose vertices are among the vertices of $P$.

For the second claim, suppose there were distinct mismatched regions
having an interior point in common.  Then their intersection should be
an interior segment for each.  Let $Q$ be one of the mismatched
regions.  It is triangulated in two different ways each coming from
the tetrahedra in one side of the hyperplane.  The triangles in either
triangulation cannot intersect improperly with the interior segment.
Thus the two triangulations of $Q$ have an interior diagonal edge in
common.  This means the triangles in $Q$ consists of more than one
connected components of the auxiliary graph, contradicting the
definition of mismatched region.
\end{proof}

\begin{proposition}
  \label{bounds}
\begin{enumerate}

\item The size of a mismatching dissection $D$ of a convex
  $3$-polytope with $n$ vertices is at least $n-2$.

\item The size of a dissection of a $3$-polytope with $n$ 
vertices is bounded from above by ${n-2 \choose 2}$.

\end{enumerate}
\end{proposition}

\begin{proof}
  $(1)$ Do an inflation of each mismatched region. This produces as
  many holes as mismatched regions, say $m$ of them. Each hole is
  bounded by two triangulations of a polygon. This is guaranteed by
  the previous lemma. Denote by $k_i$ the number of vertices of the
  polygon associated to the $i$-th mismatched region. In each of the
  holes introduce an auxiliary interior point. The point can be used
  to triangulate the interior of the holes by {\em filling in} the
  holes with the coning of the vertex with the triangles it sees.  We
  now have a triangulated ball.

  Denote by $|D|$ the size of the original dissection. The
  triangulated ball has then $|D|+ \sum_{i=1}^m 2(k_i-2)$
  tetrahedra in total. The number of interior edges of this
  triangulation is the number of interior edges in the dissection,
  denoted by $e_i(D)$, plus the new additions, for each hole of length
  $k_i$ we added $k_i$ interior edges. In a triangulation $T$ of a
  $3$-ball with $n$ boundary vertices and $n'$ interior vertices, the
  number of tetrahedra $|T|$ is related to the number of interior
  edges $e_i$ of $T$ by the formula: $|T| =n+e_i -n'-3$.  The proof is
  a simple application of Euler's formula for triangulated 2-spheres
  and 3-balls and we omit the easy details.

  Thus, we have the following equation:

$$|D|+ \sum_{i=1}^m 2(k_i-2)=n+e_i(D)+\sum_{i=1}^m k_i -m-3.$$
This can be rewritten as $|D|=n+e_i(D)-\sum_{i=1}^m k_i+3m-3$.
Taking into account that $e_i(D) \geq \sum_{i=1}^m 2(k_i-3)$
(because diagonals in a polygon are interior edges
of the dissection), we get an inequality

$$|D| \geq n +\sum_{i=1}^m k_i -3m-3.$$

Finally note that in a mismatching dissection we have $m \geq 1$ and
$k_i \geq 4$. This gives the desired lower bound.

$(2)$ Now we look at the proof of the upper bound on dissections.
Given a $3$-dissection, we add tetrahedra of volume zero to
complete to a triangulation with flat simplices that has the same
number of vertices. One can also think we are {\em filling in} the
holes created by an inflation with (deformed) tetrahedra.
  
  The lemma states that mismatched regions were of the shape of convex
  polygons.  The $2$-simplices forming a mismatched region were
  divided into two groups (those becoming apart by an inflation). The
  two groups formed different triangulations of a convex polygon,
  and they had no interior edges in common.  
  In this situation, we can make a sequence of flips (see 
  \cite{handbook}) between the two triangulations with the property
  that any edge once disappeared does not appear again (see
  Figure \ref{figure:flip}).  We add one abstract, volume zero
  tetrahedron for each flip, and obtain an abstract triangulation of a
  $3$-ball.

  \begin{figure}[htb]
\begin{center}
       \includegraphics[scale=.46]{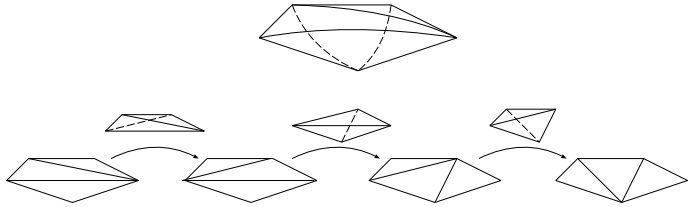}
\end{center}
\caption{Filling in holes with tetrahedra according to flips.}\label{figure:flip}
  \end{figure}
  
  The triangulation with flat simplices we created is a triangulated
  $3$-ball with $n$ vertices. By adding a new point in a fourth dimension, 
   and coning from the boundary $2$-simplices to the point, we obtain a 
triangulated $3$-sphere containing the original $3$-ball in its boundary.
 From the upper bound theorem for spheres
  (for an introduction to this topic see \cite{ziegler}) its size is
  bounded from above by the number of facets of a cyclic $4$-polytope
    minus $2n-4$, the number of $2$-simplices in the boundary of $D$.
  The 4-dimensional cyclic polytope with $n+1$ vertices is well-known to have 
  $(n+1)(n-2)/2$ facets (see \cite[page 63]{grunbaum}), 
  which completes the proof after a trivial algebraic calculation.
\end{proof}

\begin{pbm}
What is the correct upper bound  theorem for dissections
of $d$-dimensional polytopes with $d \geq 4$?
\end{pbm}

In our proof of Proposition \ref{bounds} we built a triangulated
PL-ball from a three-dimensional dissection, using the flip
connectivity of triangulations of a convex $n$-gon. Unfortunately the
same cannot be applied in higher dimensions as the flip connectivity
of triangulations of $d$-polytopes is known to be false for convex
polytopes in general \cite{santos}.  But even worse, the easy property
we used from Lemma \ref{2dmismatch} that mismatched regions are convex
polyhedra fails in dimension $d \geq 4$.

\begin{proposition} \label{badb}
 The mismatched regions of a dissection of a convex $4$-polytope 
 can be non-convex polyhedra.
\end{proposition}
\begin{proof}
The key idea is as follows: suppose we have a 3-dimensional convex polytope
$P$ and two triangulations $T_1$ and $T_2$ of it with the following
properties: removing from $P$ the tetrahedra that $T_1$ and $T_2$ have
in common, the rest is a non-convex polyhedron $P'$ such that the
triangulations $T_1'$ and $T_2'$ of it obtained from $T_1$ and $T_2$
do not have any interior 2-simplex in common (actually, something
weaker would suffice: that their common interior triangles, if any, do
not divide the interior of the polytope).

In these conditions, we can construct the dissection we want as a
bipyramid over $P$, coning $T_1$ to one of the apices and $T_2$ to the
other one. The bipyramid over the non-convex polyhedron $P'$ will be a
mismatched region of the dissection.

For a concrete example,
start with Sch\"onhardt's polyhedron whose vertices are labeled
$1,2,3$ in the lower face and $4,5,6$ in the top face.
This is a non-convex polyhedron made, for example, by twisting the three
vertices on the top of a triangular prism.
Add two
antipodal points $7$ and $8$ close to the ``top'' triangular facets
(those not breaking the quadrilaterals see Figure \ref{3dmismatch}).
For example, take as coordinates for the points
  $1=(10,0,0)$, $ 2=(-6,8,0)$, $3=(-6,-8,0)$,
  $4=(10,-0.1,10)$, $5=(-6.1,8,10)$, $6=(-5.9,-8.1,10)$,
  $7=(0,0,10.1)$,  $ 8=(0,0,-0.1).$

  \begin{figure}[htb]
    \begin{center}
        \includegraphics[scale=.4]{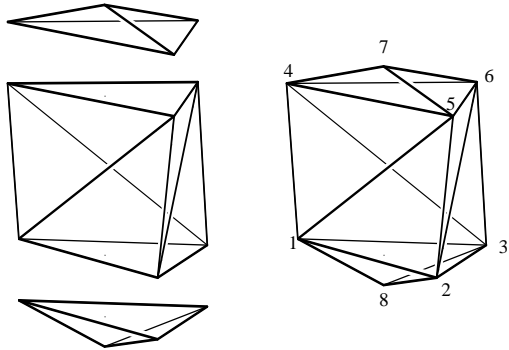}
    \end{center}
  \caption{ The mismatched region of a four-dimensional dissection.} 
 \label{3dmismatch} 
  \end{figure}

 Let $P'$ be this non-convex polyhedron and let $T'_1=
 \{1278$, $1378$, $2378$, $1247$, 
 $2457$, $2357$, $3567$, $1367$, $1467\}$ and $T'_2=\{4578$, $4678$, $5678$, 
 $1248$, $2458$, $2358$, $3568$, $1368$, $1468\}$. $T'_1$ cones vertex 7
 to the rest of the boundary of $P'$, and $T'_2$ vertex 8.
 Any common interior triangle of $T'_1$ and $T'_2$ would use the edge
 78. But the link of 78 in $T'_1$ contains only the points 1, 2 and 3, and
 the link in $T'_2$ contains only 4, 5 and 6.

 Let $P$ be the convex hull of the eight points, and let $T_1$ and
 $T_2$ be obtained from $T'_1$ and $T'_2$ by adding the three
 tetrahedra $1245$, $2356$ and $1346$. 
\end{proof}
  
\section{Optimal dissections for specific polytopes}

\begin{figure}[htb]
\begin{center}
 \includegraphics[scale=.3]{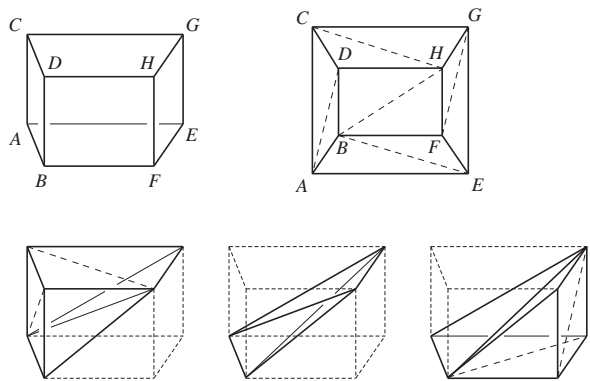} \end{center}
\caption{A triangulation of a combinatorial $3$-cube into seven tetrahedra.}
\label{cubetriangulation}
\end{figure}

The regular cube has been widely studied for its smallest dissections
\cite{haiman,hughes}. This receives the name of {\em simplexity} of
the cube. In contrast, because of the type of simplices inside a
regular $d$-cube, a simple volume argument shows that the maximal size
of a dissection is $d!$, the same as for triangulations. On the other
hand, we know that the size of the maximal triangulation of a {\em
  combinatorial} cube can be larger than that: For example, the
combinatorial $3$-cube obtained as the prism over a trapezoid
(vertices on a parabola for instance) has triangulations of size $7$.
Figure \ref{cubetriangulation} shows a triangulation with $7$ simplices
for those coordinatizations where the edges $AB$ and $GH$ are not
coplanar.  The tetrahedron $ABGH$ splits the polytope into two
non-convex parts, each of which can be triangulated with three
simplices. To see this, suppose that our polytope is a very small
perturbation of a regular $3$-cube. In the regular cube, $ABGH$ becomes
a diagonal plane which divides the cube into two triangular prisms
$ABCDGH$ and $ABEFGH$.  In the non-regular cube, the diagonals $AH$
and $BG$, respectively, become non-convex. Any pair of triangulations
of the two prisms, each using the corresponding diagonal,
together with tetrahedron $ABGH$ give a
triangulation of the perturbed cube with $7$ tetrahedra.  The boundary
triangulation is shown in the flat diagram.  It is worth noticing that
for the regular cube the boundary triangulation we showed does not
extend to a triangulation of the interior.

  One can then ask, what is the general growth for the size of a
  maximal dissection of a combinatorial cube? To answer this question,
  at least partially, we use the above construction and we adapt an
  idea of M. Haiman, originally devised to produce 
  small triangulations of regular cubes
  \cite{haiman}.  The idea is that 
  from triangulations of a $d_1$-cube and a $d_2$-cube
  of sizes $s_1$ and $s_2$ respectively we can get 
  triangulations of the $({d_1+d_2})$-cube by first subdividing it into $s_1
  \times s_2$ copies of the product of two simplices of dimensions
  $d_1$ and $d_2$ and then 
  triangulating each such piece. We recall that any triangulation of the
  Cartesian product of a $d_1$-simplex and a $d_2$-simplex has
  ${d_1+d_2 \choose d_1}$ maximal simplices.  Hence, in
  total we have a triangulation of the $({d_1+d_2})$-cube into $s_1 \times s_2
  \times {d_1+d_2 \choose d_1}$ maximal simplices.  
  Recursively, if one starts with a triangulation of size $s$ of the $d$-cube,
  one obtains triangulations for the $rd$-cube of size 
  $(rd)! ({s \over d!})^r$. In Haiman's context one wants $s$ to be
  small, but here we want it to be big.

  More precisely, denote by $f(d)$ the function $\max_{C:\
  d\mbox{\scriptsize-cube}}(\max_{T\mbox{\scriptsize\ of \,}C\,}
  |T|)$ and call $g(d)=\left(f(d)/d!\right)^{1/d}$.
  Haiman's argument shows that if $f(d_1) \ge
  {c_1}^{d_1} {d_1}!$ and $f(d_2) \ge {c_2}^{d_2} {d_2}!$ for
  certain constants $c_1$ and $c_2$ then 
  $f(d_1 +   d_2) \ge {c_1}^{d_1}{c_2}^{d_2} ({d_1+d_2})!$.
  Put differently, 
  that $g(d_1+d_2)\ge \left(g(d_1)^{d_1} g(d_2)^{d_2}\right)^{1/(d_1+d_2)}$.
  The value on the right hand side is the weighted geometric mean of
  $g(d_1)$ and $g(d_2)$. In particular, if both $g(d_1)$ and
  $g(d_2)$ are $\ge 1$
  and one of them is $>1$ then $g(d_1+d_2)$ is $>1$ as well.

  We have constructed above a triangulation of size 7
  for the Klee-Minty $3$-cube, which proves $g(3)\ge
  \sqrt[3]{7/6}=1.053$. With Haiman's idea 
  we can now construct ``large'' triangulations of certain $4$-cubes
  and $5$-cubes, which prove respectively that $g(4)\ge \sqrt[4]{7/6}=1.039$
  and $g(5)\ge \sqrt[5]{7/6}=1.031$ 
  (take $d_1=3$ and $d_2$ equal to one and two
  respectively). Finally,
  since any $d>5$ can be expressed as a sum of 3's and 4's, we have
  $g(d)\ge\min\{g(3),g(4)\}\ge 1.039$ for any $d>5$. Hence:

\begin{proposition} \label{combinatorialcube}
For the family of combinatorial $d$-cubes with $d>2$ the function 
{\rm $ f(d)= \max_{C:\
d\mbox{\scriptsize-cube}}(\max_{T\mbox{\scriptsize\ of \,}C}\, |T|)$ 
}
admits the lower bound $f(d) \geq c^d d!$ where $c\ge 1.031$.
\end{proposition}

Exactly as in Haiman's paper, the constant $c$ can be improved
(asymptotically) if one
starts with larger triangulations for the smaller dimensional cubes. 
Using computer
calculations (see Remark \ref{rem.dhss}), we obtained a maximal
triangulation for the Klee-Minty 4-cube with 38 maximal simplices,
which shows that $g(d)\geq\sqrt[4]{38/24}=1.122$ for every $d$
divisible by 4 (see \cite{amezie} for a
complete study of this family of cubes). We omit listing the triangulation
here but it is available from the authors by request.


\begin{pbm}
Is the sequence $g(d)$ bounded? In other words, 
is there an upper bound of type $c^d d!$ for the function $f(d)$?
Observe that the same question for {\em minimal} triangulations of the
{\em regular} $d$-cube $($whether there is a lower bound of type $c^d
d!$ for some $c>0$$)$ is open as well. See \cite{smith} for the best
lower bound known.
\end{pbm}

We continue our discussion with the study of optimal triangulations
for three-dimensional prisms and antiprisms. 
We will call an {\em $m$-prism} any 3-polytope with the combinatorial type
of the product of a convex $m$-gon with a line segment.  
An {\em $m$-antiprism} will be any 3-polytope whose faces are two
convex $m$-gons and $2m$ triangles, each $m$-gon being
adjacent to half of the triangles.
Vertices of the two $m$-gons are connected with a band of
alternately up and down pointing triangles. 

Each such polyhedron has a regular coordinatization in which all the faces are
regular polygons, and a
realization space which is the set of all possible coordinatizations
that yield the same combinatorial information \cite{juergen}. Our
first result is valid in the whole realization space.

\begin{proposition} \label{prismantip}
  For any three-dimensional $m$-prism, in any of its possible
  coordinatizations, the number of tetrahedra in a minimal
  triangulation is $2m-5+\lceil {m\over 2}\rceil$. 

  For any three-dimensional $m$-antiprism, in any of its
  possible coordinatizations, the number of tetrahedra in a minimal
  triangulation is $3m-5$.
\end{proposition}
\begin{proof} 
In what follows we use the word {\em cap} to refer to the $m$-gon
facets appearing in a prism or antiprism. We begin our discussion
proving that any triangulation of the prism or antiprism has at least
the size we state, and then we will construct triangulations with
exactly that size.

We first prove that every triangulation of the $m$-prism requires at
least $2m-5+\lceil {m\over 2}\rceil$ tetrahedra. We call a tetrahedron
of the $m$-prism {\em mixed} if it has two vertices on the top cap and
two vertices on the bottom cap of the prism, otherwise we say that the
tetrahedron is {\em top-supported} when it has three vertices on the top
(respectively {\em bottom-supported}). For example, Figure
\ref{fig.prism} shows a triangulation of the regular 12-prism, in three
slices. Parts (a) and (c) represent, respectively, the bottom and top
caps. Part (b) is the intersection of the prism with the parallel
plane at equal distance to both caps. In this intermediate slice,
bottom or top supported tetrahedra appear as triangles, while mixed
tetrahedra appear as quadrilaterals.

\begin{figure}[htb]
  \begin{center}
        \includegraphics[scale=.4]{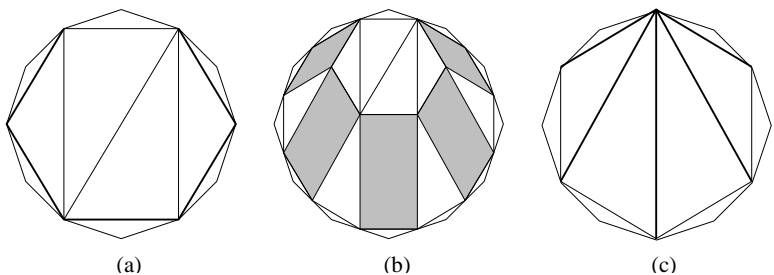}
  \end{center}
\caption{A minimal triangulation of the regular 12-prism.}
 \label{fig.prism}
\end{figure}

Because all triangulations of
an $m$-gon have $m-2$ triangles there are always exactly $2m-4$ tetrahedra
that are bottom or top supported.  In the rest, we show there are at
least $\lceil{m\over 2}\rceil-1$ mixed tetrahedra.  Each mixed
tetrahedra {\em marks} an edge of the top, namely the edge it uses
from the top cap. Of course, several mixed tetrahedra could mark the
same top edge. Group together top-supported tetrahedra that have the
same bottom vertex.
This grouping breaks the triangulated top $m$-gon
into polygonal regions. Note that every edge between two of these
regions must be marked. For
example, in part (c) of Figure \ref{fig.prism} the top cap is divided
into 6 regions by 5 marked edges (the thick edges
in the Figure). Let $r$ equal the number
of regions under the equivalence relation we set. There are $r-1$
interior edges separating the $r$ regions, and all of them are marked. Some
boundary edges of the top cap may be marked too (none of them is
marked in the example of Figure \ref{fig.prism}).

We can estimate the marked edges in another way: There are $m$ edges
on the boundary of the top, which appear partitioned among some of the
regions (it could be the case some region does not contain any
boundary edge of the $m$-gon). We claim that no more than {\em two}
boundary edges
per region will be unmarked $(\ast)$. This follows because a boundary
edge is not marked only when the top supported tetrahedron that
contains it has the point in the bottom cap that is directly under
one of the vertices of the edge. In a region, at most two boundary edges can
satisfy this.  Hence we get at least $m-2r$ marked edges on the
boundary of the top and at least $(r-1)+(m-2r)=m-r-1$ marked edges
in total. Thus the number of mixed tetrahedra is at least the maximum of
$r-1$ and $m-r-1$. In conclusion, we get that, indeed, the number of
mixed tetrahedra is bounded below by $\lceil {m\over 2}\rceil-1$. Note
that we only use the combinatorics and convexity of the prism in our
arguments.
We will show that minimal triangulations achieve this lower bound,
but then, observe that if $m$ is even, 
in a minimal triangulation we must have
$r=m/2$ and no boundary edge can be marked, as is the case
in Figure \ref{fig.prism}. If $m$ is odd, then we must have 
$r\in \{(m-1)/2, (m+1)/2\}$
and  at most one boundary edge can be marked.

The proof that any triangulation of an $m$-antiprism includes at least
$3m-5$ tetrahedra is similar.  There are
$2m-4$ top-supported and bottom-supported tetrahedra in any
triangulation and there are $r-1$ marked edges between the regions in
the top.  The only difference is that, instead of claim $(\ast)$, one
has at most {\em one} unmarked boundary edge per region.  Thus there are
at least $m-r$ marked edges in the boundary of the top, and in total
at least $(r-1)+(m-r)=m-1$ marked edges in the top.  Hence there exist
at least $(2m-4)+(m-1)=3m-5$ tetrahedra in any triangulation.

For an $m$-antiprism we can easily create a triangulation of size
$3m-5$ by choosing any triangulation of the bottom $m$-gon and then
coning a chosen vertex $v$ of the top $m$-gon to the $m-2$ triangles
in that triangulation and to the $2m-3$ triangular facets of the
$m$-antiprism which do not contain $v$. This construction is exhibited
in Figure \ref{fig.antiprism}. Parts (a) and (c) show the bottom and
top caps triangulated (each with its $5$ marked edges) and part (b)
an intermediate slice with the 5 mixed tetrahedra appearing as quadrilaterals.

\begin{figure}[htb]
  \begin{center}
        \includegraphics[scale=.4]{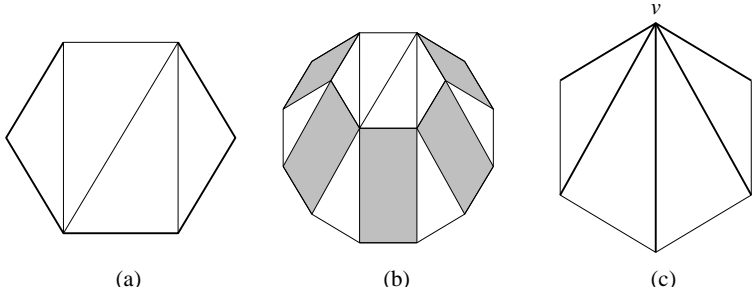}
  \end{center}
\caption{A minimal triangulation of the regular 6-antiprism.}
 \label{fig.antiprism}
\end{figure}

For an $m$-prism, let $u_i$ and $v_i$, $i=1,\dots,m$ denote the top
and bottom vertices respectively, so that the vertices of each cap are
labeled consecutively and $u_iv_i$ is always an edge of the prism.

If $m$ is even we can chop off the vertices $u_i$ for odd $i$ and
$v_j$ for even $j$, so that the prism is decomposed into $m$
tetrahedra and an $({m\over 2})$-antiprism. The antiprism can be
triangulated into ${3m\over 2}-5$ tetrahedra, which gives a
triangulation of the prism into ${5m\over 2}-5$ tetrahedra, as
desired. Actually, this is how the triangulation of Figure
\ref{fig.prism} can be obtained from that of Figure \ref{fig.antiprism}.

If $m$ is odd we do the same, except that we chop off only the
vertices $u_1,\dots,u_{m-2}$ and $v_2,\dots,v_{m-1}$ (no vertex is
chopped in the edge $u_mv_m$). This produces $m-1$ tetrahedra and an
$({m+1\over 2})$-antiprism. We triangulate the antiprism into
${3m+3\over 2}-5$ tetrahedra and this gives a triangulation of the
$m$-prism into ${5m+1\over 2}-5$ tetrahedra.
 \end{proof}

We have seen that the coordinates are not important when calculating
minimal triangulations of the three-dimensional prisms and
antiprisms. On the other hand, the difference in size of the maximal
triangulation can be quite dramatic. Below we prove that in 
certain coordinatizations it is roughly ${m^2\over 2}$ and show
experimental data indicating that for the regular prism it is close to
${m^2\over 4}$.

\begin{proposition}
\label{maxprism}
Let $A_m$ be a prism of order $m$, with all its side edges parallel.
\begin{enumerate}
\item The size of a maximal triangulation of $A_m$ is bounded as
$$ \left\lceil\frac{m^2+6m-16}{4}\right\rceil
\leq \max_{T\mbox{\scriptsize\ of }A_m}|T| \leq \frac{m^2+m-6}{2}.$$

\item The upper bound is achieved if the two caps 
$($$m$-gon facets$)$ are parallel and
there is a direction in which the whole prism projects onto one of its
side quadrangular facets. $($For a concrete example, let one of the
$m$-gon facets have vertices on a parabola and let $A_m$ be the
product of it with a segment$)$.
\end{enumerate}
\end{proposition}

\begin{proof}
  Let the vertices of the prism be labeled $u_1,\dots,u_m$ and
  $v_1,\dots,v_m$ so that the $u_i$'s and the $v_j$'s form the two
  caps, vertices in each cap are labeled consecutively and $u_iv_i$
  is always a side edge.

  For the upper bound in part (1), we have to prove that a
  triangulation of $A_m$ has at most
  $\frac{m^2+m-6}{2}-2m+3={m(m-3)\over 2}$ interior diagonals. The
  possible diagonals are the edges $u_iv_j$ where $i-j$ is not in
  $\{-1,0,1\}$ modulo $m$. This gives exactly twice the number we
  want. But for any $i$ and $j$ the diagonals $u_iv_j$ and $u_jv_i$
  intersect, so only one of them can appear in each triangulation.

  We now prove that the upper bound is achieved if $A_m$ is in the
  conditions of part (2). In fact, the condition on $A_m$ that we will
  need is that for any $1\le i<j\le k<l\le m$, the point $v_j$ sees
  the triangle $v_iu_ku_l$ from the same side as $v_k$ and $v_l$
  (i.e.\ ``from above'' if we call top cap the one containing the
  $v_i$'s). With this we can construct a triangulation with
  $\frac{m^2+m-6}{2}={m-1\choose 2}+2m-4$ tetrahedra, as follows:

  First cone the vertex $v_1$ to any triangulation of the bottom cap
  (this gives $m-2$ tetrahedra). The $m-2$ upper boundary facets of
  this cone are visible from $v_2$, and we cone them to it (again
  $m-2$ tetrahedra). The new $m-2$ upper facets are visible from $v_3$
  and we cone them to it ($m-2$ tetrahedra more). Now, one of the
  upper facets of the triangulation is $v_1v_2v_3$, part of the upper
  cap, but the other $m-3$ are visible from $v_4$, so we cone them and
  introduce $m-4$ tetrahedra. Continuing the process, we will
  introduce $m-4$, $m-5,\dots,2,1$ tetrahedra when coning the vertices
  $v_5,v_6,\dots,v_{m-1},v_m$, which gives a total of ${m-1\choose
    2}+2m-4$ tetrahedra, as desired.

  The triangulation we have constructed is the {\em placing
    triangulation} \cite{handbook} associated to any ordering of the
  vertices finishing with $v_1,\dots,v_m$. A different description of
  the same triangulation is that it cones the bottom cap to $v_1$, the
  top cap to $u_m$, and its mixed tetrahedra are all the possible
  $v_{i}v_{i+1}u_ju_{j+1}$ for $1\le i < j\le m-1$. This gives
  ${m-1\choose 2}$ mixed tetrahedra, and ${m-1\choose 2}+2m-4$
  tetrahedra in total.

  We finally prove the lower bound stated in part (1). Without loss of
  generality, we can assume that our prism has its two caps parallel
  (if not, do a projective transformation keeping the side edges
  parallel). Then, $A_m$ can be divided into two prisms in the
  conditions of part (2) of sizes $k$ and $l$ with $k+l=m+2$: take any
  two side edges of $A_m$ which posses parallel supporting planes and
  cut $A_m$ along the plane containing both edges.  By part (2), we
  can triangulate the two subprisms with ${k+1 \choose 2}-3$ and
  ${l+1\choose 2}-3$ tetrahedra respectively, taking care that the two
  triangulations use the same diagonal in the dividing plane. This
  gives a triangulation of $A_m$ with ${k+1 \choose 2}+{l+1\choose
    2}-6= {k^2+l^2+m-10\over 2}$ tetrahedra.  This expression achieves
  its minimum when $k$ and $l$ are as similar as possible, i.e.\ 
  $k=\lfloor {m\over 2}\rfloor +1$ and $l=\lceil {m\over 2}\rceil +1$.
  Plugging these values in the expression gives a triangulation of
  size $\left\lceil\frac{m^2+6m-16}{4}\right\rceil$.
\end{proof}

Based on an integer programming approach we can compute maximal triangulations
of specific polytopes (see remark at the end of the article).
Our computations with regular prisms up to $m=12$ show that the size
of their maximal triangulations achieve the lower bound stated in part
(1) of Proposition \ref{maxprism} (see Table \ref{tableprisms}).  In
other words, that the procedure of dividing them into two prisms of
sizes $\lfloor {m\over 2}\rfloor +1$ and $\lceil {m\over 2}\rceil +1$
in the conditions of part (2) of Proposition \ref{maxprism} and
triangulating the subprisms independently yields maximal
triangulations. 

We have also computed maximal sizes of triangulations for the
regular $m$-antiprisms up to $m=12$, which turn out to follow the formula
$\left\lfloor\frac{m^2+8m-16}{4}\right\rfloor$.
A construction of a triangulation of this size for every $m$
can be made as follows:
Let the vertices of the regular $m$-antiprism be labeled
$u_1,\dots,u_m$ and $v_1,\dots,v_m$ so they are forming the
vertices of the two caps consecutively in this order and $v_iu_i$ and
$u_iv_{i+1}$ are side edges. We let $v_{m+1}=v_1$. The triangulation
is made by placing the vertices in any ordering finishing with
$v_1,v_2,v_m,v_3,v_{m-1},\dots,v_{\lceil\frac{m}{2}\rceil+1}$. The
tetrahedra used are the bottom-supported tetrahedra with
apex $v_1$, top-supported tetrahedra with apex
$u_{\lceil\frac{m}{2}\rceil}$ and the mixed tetrahedra
$v_{i}v_{i+1}u_ju_{j+1}$ for $1\le i \le j \le \lfloor\frac{m}{2}\rfloor$ and
$u_iu_{i+1}v_jv_{j+1}$ for $\lfloor\frac{m}{2}\rfloor+1\leq i<j\leq
m$.

We conjecture that these formulas for regular base prisms and
antiprisms actually give the sizes of their maximal triangulations for
every $m$,
but we do not have a proof.

\begin{table}[htpb]
\begin{center}
\small
\font\eightrm=cmr8
\eightrm
\begin{tabular}{|c|c|c|c|c|c|c|c|c|c|c|}
\hline
         $m$             & 3 & 4 &  5 &  6 &  7 &  8 &  9 & 10 & 11 & 12 \cr
\hline
Prism (regular base)     & 3 & 6 & 10 & 14 & 19 & 24 & 30 & 36 & 43 & 50 \cr
\hline
Antiprism (regular base) & 4 & 8 & 12 & 17 & 22 & 28 & 34 & 41 & 48 & 56 \cr
\hline
\end{tabular}
\end{center}
\caption{Sizes of maximal triangulations of prisms and antiprisms.}
\label{tableprisms}
\end{table}

\begin{rem}
\label{rem.dhss}
\rm How can one find minimal and maximal triangulations in specific
instances? The approach we followed for computing Tables
\ref{res_table} and \ref{tableprisms} and some of the results in
Proposition \ref{prop.latticepoly} is the one proposed in \cite{dhss},
based on the solution of an integer programming problem. We think of
the triangulations of a polytope as the vertices of the following
high-dimensional polytope: Let $A$ be a $d$-dimensional polytope with
$n$ vertices.  Let $N$ be the number of $d$-simplices in $A$.  We
define $P_A$ as the convex hull in $\bbbr^N$ of the set of incidence
vectors of all triangulations of ${A}$. For a triangulation $T$ the
{\em incidence vector} $v_T$ has coordinates $\,(v_T)_\sigma = 1$ if
$\sigma \in T$ and $\,(v_T)_\sigma = 0$ if $\sigma \not\in T$.  The
polytope $P_A$ is the {\em universal polytope} defined in general by
Billera, Filliman and Sturmfels \cite{BIFIST} although it appeared in
the case of polygons in \cite{dantzig}. In \cite{dhss}, it was shown
that the vertices of $P_A$ are precisely the integral points inside a
polyhedron that has a simple description in terms of the oriented
matroid of $A$ (see \cite{dhss} for information on oriented
matroids). The concrete integer programming problems were solved using
{\em C-plex Linear Solver}$^{TM}$.  The program to generate the linear
constraints is a small $C^{++}$ program written by Samuel Peterson and
the first author. Source code, brief instructions, and data files are
available via ftp at {\bf http://www.math.ucdavis.edu/\~{}deloera}. An
alternative implementation by A. Tajima is also available
\cite{tajima1,tajima2}. He used his program to corroborate some of
these results.

It should be mentioned that a simple variation of the ideas in
\cite{dhss} provides enough equations for an integer program whose
feasible vertices are precisely the $0/1$-vectors of dissections.  The
incidence vectors of dissections of $conv(A)$, for a point set $A$,
are just the $0/1$ solutions to the system of equations $\langle x,
v_T \rangle = 1$, where $v_T$'s are the incidence vectors for every
regular triangulation $T$ of the Gale transform $A^*$ (regular
triangulations in the Gale transform are the same as chambers in $A$).
Generating all these equations is as hard as enumerating all the
chambers of $A$. Nevertheless, it is enough to use those equations
coming from placing triangulations (see \cite[Section
3.2]{triangOMS}), which gives a total of about $n^{d+1}$ equations if
$A$ has $n$ points and dimension $d$.
\end{rem}

\section*{Acknowledgments}
We are grateful to Alexander Below and J\"urgen Richter-Gebert for
their help and ideas in the proofs of Proposition \ref{bounds} and
\ref{badb}. Alexander Below made Figure \ref{3dmismatch} using the
package Cinderella. The authors thank Akira Tajima and J\"org Rambau
for corroborating many of the computational results. We thank Samuel
Peterson for his help with our calculations. Finally, we thank Hiroshi
Imai, Bernd Sturmfels, and Akira Tajima for their support of this
project.

\end{document}